\documentclass[preprint,11pt]{elsarticle}
\usepackage[T1]{fontenc}
\usepackage{lmodern,mathtools,amssymb,amsthm}
\usepackage[hidelinks]{hyperref}
\usepackage{microtype}
\hypersetup{pdftitle={Parity Balancing of ℓ₁ Unit Vectors},pdfauthor={Chunyin Zhao, Muwen Li}}

\newtheorem{theorem}{Theorem}[section]
\newtheorem{lemma}[theorem]{Lemma}
\newtheorem{proposition}[theorem]{Proposition}

\theoremstyle{definition}

\newtheorem{remark}[theorem]{Remark}

\newcommand{\R}{\mathbb{R}}
\newcommand{\eps}{\varepsilon}
\newcommand{\one}{\mathbf{1}}
\newcommand{\sgn}{\operatorname{sgn}}
\newcommand{\rank}{\operatorname{rank}}

\allowdisplaybreaks

\begin{document}
\begin{frontmatter}
\title{Parity Balancing of $\ell_1$ Unit Vectors}
\author[btbu]{Chunyin Zhao}
\author[ustc]{Muwen Li}
\address[btbu]{Beijing Technology and Business University, Liangxiang Higher Education Park, Fangshan District, Beijing 102488, P.R. China}
\address[ustc]{School of Mathematical Sciences, University of Science and Technology of China, No.~96 Jinzhai Road, Hefei, Anhui 230026, P.R. China}

\begin{abstract}
Let \(v_1,\dots,v_n\) be unit vectors in \(\ell_1^d\), and suppose that \(n\) and
\(d\) have opposite parity. We prove that for every \(d\ge 2\) there are signs
\(\eps_1,\dots,\eps_n\in\{-1,1\}\) such that
\(\|\sum_{i=1}^n \eps_i v_i\|_1\le d-1\), and that \(d-1\) is sharp.
The proof refines the classical rounding argument. A
local rigidity lemma shows that if every corner of a certain parallelepiped lies
outside the \(\ell_1\)-ball of radius \(d-1\), then the generating vectors have
pairwise distinct dominant coordinates. Adjacent vertices exchange one free
vector while preserving the number of positive fixed signs; connectedness and
central symmetry give the parity contradiction. We also give deterministic
constructions: a direct one in \(O(nd^2)\) exact arithmetic operations and
comparisons, and a blocked one in
\(O(nd^{(\omega+1)/2}+d^\omega)\), where \(\omega\) is the matrix multiplication
exponent. With \(\omega\le9/4\), the latter is \(O(nd^{13/8})\) when \(n>d\).
\end{abstract}

\begin{keyword}
\(\ell_1\)-norm \sep balancing unit vectors \sep parity \sep deterministic algorithm
\MSC[2020] 05D40 \sep 52C07 \sep 68W30
\end{keyword}
\end{frontmatter}

\section{Introduction}

Given unit vectors in a normed space, one seeks signs whose sum is close to the
origin. In \(\ell_1^d\), the classical rounding argument gives \(d\): relax the
signs to \([-1,1]\), take a vertex of the resulting polytope, and round its free
coordinates. The coordinate basis shows that \(d\) is best possible.

Repeating a basis vector changes the parity: two equal vectors can cancel.
Swanepoel proved that an odd number of unit vectors in any normed plane admits
signs whose sum has norm at most \(1\), and asked for the higher-dimensional
analogue \cite{Swanepoel}; the case \(d=2\) below is a special case. We prove
the following parity improvement (see also \cite{Spencer}).

\begin{theorem}\label{thm:main}
Let \(d\ge 2\), let \(n\not\equiv d\pmod 2\), and let
\(v_1,\dots,v_n\in\R^d\) satisfy \(\|v_i\|_1=1\). Then there exist
\(\eps_1,\dots,\eps_n\in\{-1,1\}\) such that
\begin{equation}\label{eq:main}
\left\|\sum_{i=1}^n \eps_i v_i\right\|_1\le d-1.
\end{equation}
\end{theorem}

The bound \(d-1\) is sharp: for \(e_1,\dots,e_{d-1},e_d,e_d\), the first
\(d-1\) coordinates have absolute value \(1\) in every signed sum. The
unit-vector assumption and the choice of target \(0\) are also essential: a zero
vector would make \(e_1,\dots,e_d,0\) violate the parity conclusion, and in
dimension \(3\) the vectors \(e_1,e_1,e_2,e_2,e_3,e_3\) have \(\ell_1\)-distance
at least \(3>d-1=2\) from \((1,1,1)\). For \(n\ge d\), the parity hypothesis is
necessary as well: the standard basis, padded with copies of \(e_1\), gives
lower bound \(d\) when \(n\equiv d\pmod2\).

The proof is by contradiction. If all sign choices fail, the rounding estimate
forces every vertex to have exactly \(d\) free coordinates; local rigidity gives
them distinct dominant coordinates. Adjacent vertices exchange one free vector,
with a positive coefficient by Cramer's rule. Hence the number of positive fixed
signs is constant on the vertex graph. Central symmetry gives the complementary
count at \(-x\), a contradiction.

Section~\ref{sec:rounding} gives vertex rounding, Section~\ref{sec:rigidity}
local rigidity, Section~\ref{sec:exchange} the exchange and parity
contradiction, and Section~\ref{sec:algorithm} the deterministic constructions.
\ref{app:height} gives a height interpretation.

\section{Vertex rounding}\label{sec:rounding}

Write \(\ell_1^d\) for \(\R^d\) equipped with the norm
\(\|x\|_1=\sum_{j=1}^d |x_j|\). For vectors \(v_1,\dots,v_n\in\R^d\), put
\(V=[v_1,\dots,v_n]\), so that
\[
Vx=\sum_{i=1}^n x_i v_i.
\]
The zero fiber of \(V\) is
\[
K(V)=\{x\in[-1,1]^n:Vx=0\}.
\]
For \(x\in K(V)\), the index \(i\) is \emph{free} if \(|x_i|<1\), and
\emph{fixed} if \(|x_i|=1\). We write
\[
F(x)=\{i:|x_i|<1\},\qquad Z(x)=\{1,\dots,n\}\setminus F(x).
\]
A point of a convex set is a \emph{vertex} if it is not the midpoint of two
distinct points of that set.

\begin{lemma}[Vertices and rounding]\label{lem:vertex-criterion}
Let \(c\in\R^d\), let \(K(V,c)=\{x\in[-1,1]^n:Vx=c\}\), and let
\(x\in K(V,c)\). Then \(x\) is a vertex of \(K(V,c)\) if and only if the
vectors \(\{v_i:i\in F(x)\}\) are linearly independent. In that case
\(|F(x)|\le d\). Moreover, if one keeps every fixed coordinate at its sign and
rounds every free coordinate to the nearest sign, obtaining a sign vector
\(\eps\), then
\begin{equation}\label{eq:rounding}
\|V\eps-c\|_1\le \sum_{i\in F(x)}(1-|x_i|)
=|F(x)|-\|x_{F(x)}\|_1.
\end{equation}
\end{lemma}

\begin{proof}
Let \(x=(y+z)/2\) with \(y,z\in K(V,c)\). At every bound coordinate
\(i\notin F(x)\), convexity forces \(y_i=z_i=x_i\). Hence
\(V_{F(x)}(y_{F(x)}-x_{F(x)})=0\).
If the free columns are dependent, choose a nonzero \(g\) supported on
\(F(x)\) with \(Vg=0\); then \(x\pm\eta g\in K(V,c)\) for small \(\eta>0\),
so \(x\) is not a vertex. If they are independent, the last equation forces
\(y_{F(x)}=x_{F(x)}\), and hence \(y=z=x\). Thus \(x\) is a vertex exactly
when the free columns are independent.

The free vectors lie in \(\R^d\), so their number is at most \(d\). For the
rounding estimate, set \(\eps_i=x_i\) for fixed \(i\), and choose
\(\eps_i\in\{-1,1\}\) nearest to \(x_i\) for free \(i\). Then
\(|\eps_i-x_i|=1-|x_i|\) on free coordinates and \(0\) on fixed coordinates.
Using \(Vx=c\),
\[
\begin{aligned}
\|V\eps-c\|_1
&=\|V(\eps-x)\|_1
=\left\|\sum_{i\in F(x)}(\eps_i-x_i)v_i\right\|_1\\
&\le\sum_{i\in F(x)}|\eps_i-x_i|\,\|v_i\|_1
=\sum_{i\in F(x)}(1-|x_i|).
\end{aligned}
\]
This proves \eqref{eq:rounding}.
\end{proof}

Since \(K(V,c)\) is a compact convex subset of \(\R^n\), the
Krein--Milman theorem \cite{KreinMilman} shows that every nonempty \(K(V,c)\)
has a vertex.

Taking a vertex \(x\) of \(K(V)\) and rounding as in
Lemma~\ref{lem:vertex-criterion} gives the classical bound
\(\|V\eps\|_1\le |F(x)|\le d\).

\section{Local rigidity and cone counting}\label{sec:rigidity}

We prove Theorem~\ref{thm:main} by contradiction. Assume throughout this and
the next section that
\begin{equation}\label{eq:H}
\|V\eps\|_1>d-1\qquad\text{for every }\eps\in\{-1,1\}^n.
\end{equation}

If \(n<d\), then \(\|V\eps\|_1\le\sum_i\|v_i\|_1=n\le d-1\), contradicting
\eqref{eq:H}. The case \(n=d\) is excluded by the parity assumption. Thus
\(n>d\).

Now let \(x\) be any vertex of \(K(V)\). By
Lemma~\ref{lem:vertex-criterion}, the number \(f\) of free coordinates is at
most \(d\). If \(f\le d-1\), rounding as in \eqref{eq:rounding} gives
\(\|V\eps\|_1\le f\le d-1\), contrary to \eqref{eq:H}. Hence \(f=d\). Write
\(F(x)=\{i_1<\dots<i_d\}\), put \(a_j=v_{i_j}\),
\(A=[a_1,\dots,a_d]\), \(y_j=x_{i_j}\), and let
\[
p=\sum_{i\notin F(x)}x_i v_i
\]
be the signed sum of the fixed vectors. The \(d\) free columns of \(A\) are
independent, so \(A\) is invertible; also \(n-d\ge1\).

The rounding estimate gives \(\|V\eps\|_1\le d-\|y\|_1\). Since the left side
is larger than \(d-1\), we get \(\|y\|_1<1\). Finally, \(Vx=0\) means
\(Ay+p=0\). If the fixed coordinates retain their signs \(x_i\) and the free
coordinates are replaced by \(\xi\), the resulting signed sum is
\[
A\xi+p=A(\xi-y).
\]
By \eqref{eq:H}, its norm is larger than \(d-1\), that is,
\begin{equation}\label{eq:all-corners-fail}
\|A(\xi-y)\|_1>d-1\qquad\text{for every }\xi\in\{-1,1\}^d.
\end{equation}
We shall use these facts repeatedly: under \eqref{eq:H}, every vertex of
\(K(V)\) has exactly \(d\) free coordinates, its free coefficient vector has
\(\ell_1\)-norm smaller than \(1\), and all of its corresponding corners fail
\eqref{eq:all-corners-fail}.

For the rigidity step, let \(a\in\R^d\) with \(\|a\|_1=1\). A coordinate \(j\)
is \emph{dominant} if \(|a_j|>1/2\). Such a coordinate is unique, because two
coordinates with absolute value larger than \(1/2\) would already give
\(\|a\|_1>1\). Associated with coordinate \(j\) and sign \(s\in\{-1,1\}\) are
the open convex cones
\begin{equation}\label{eq:cone}
\mathcal C_{j,s}=\left\{z\in\R^d: s z_j>\sum_{k\ne j}|z_k|\right\}.
\end{equation}
Thus \(j\) is a dominant coordinate of \(a\) precisely when \(a\in
\mathcal C_{j,s}\) for some \(s\).

Two elementary facts will be used. If \(z\in\mathcal C_{j,s}\), then
\begin{equation}\label{eq:cone-cross}
s z_j+t z_k>0
\end{equation}
for every \(k\ne j\) and every \(t\in\{-1,1\}\). Indeed,
\(s z_j>\sum_{l\ne j}|z_l|\ge |z_k|\ge -t z_k\). Also, each
\(\mathcal C_{j,s}\) is convex, and cones with different labels are disjoint:
if \(z\in\mathcal C_{j,s}\cap\mathcal C_{k,t}\) with \(j\ne k\), then
\(|z_j|\ge s z_j>|z_k|\) and \(|z_k|\ge t z_k>|z_j|\), a contradiction; the
case \(j=k\), \(s\ne t\) is immediate.

\begin{lemma}[Cone counting]\label{lem:cone-counting}
Let \(a_1,\dots,a_d,p\in\R^d\). Suppose that each cone \(\mathcal C_{r,s}\)
contains at least one of the \(2d\) labeled points \(p+a_i,p-a_i\), counted
with multiplicity. Then the \(a_i\) have pairwise distinct dominant positions.
\end{lemma}

\begin{proof}
Count with multiplicity. Each cone contains at least one point and each point
lies in at most one cone; since there are \(2d\) of each, the correspondence is
a bijection.

For each \(i\), let
\[
p+a_i\in\mathcal C_{r_i,s_i},\qquad p-a_i\in\mathcal C_{q_i,t_i}.
\]
Define \(D_i=s_i p_{r_i}+t_i p_{q_i}\). If \(r_i=q_i\), the two cones are
different, so \(t_i=-s_i\) and \(D_i=0\). If \(r_i\ne q_i\), apply
\eqref{eq:cone-cross} to \(p+a_i\) at coordinate \(q_i\) with sign \(t_i\), and
to \(p-a_i\) at coordinate \(r_i\) with sign \(s_i\). Adding the two strict
inequalities cancels the terms containing \(a_i\) and gives \(2D_i>0\). Thus
always \(D_i\ge0\), with strict inequality when \(r_i\ne q_i\).

The sum \(\sum_iD_i\) assigns to a point in \(\mathcal C_{r,s}\) the number
\(s p_r\). By the bijection, it equals
\[
\sum_{r=1}^d\left(p_r+(-p_r)\right)=0.
\]
Hence every \(D_i=0\), so \(r_i=q_i\) and \(t_i=-s_i\). It follows that
\(p+a_i\in\mathcal C_{r_i,s_i}\) and \(a_i-p=-(p-a_i)\in\mathcal C_{r_i,s_i}\).
By convexity, their midpoint \(a_i\) lies in \(\mathcal C_{r_i,s_i}\). Thus
each \(a_i\) has dominant position \(r_i\), and different \(i\) use different
pairs of opposite cones, so the positions are pairwise distinct.
\end{proof}

\begin{lemma}[Local rigidity]\label{lem:rigidity}
Let \(A=[a_1,\dots,a_d]\) with \(\|a_i\|_1=1\), and let
\(y\in[-1,1]^d\). Suppose that
\[
\|A(\xi-y)\|_1>d-1\qquad\text{for every }\xi\in\{-1,1\}^d.
\]
Then every \(a_i\) has a dominant coordinate, and these dominant coordinates
have pairwise distinct positions.
\end{lemma}

\begin{proof}
We first reduce to the case in which \(A\) is invertible and
\(y\in(-1,1)^d\). If \(A\) is singular, the nonempty compact convex fiber
\(K'=\{z\in[-1,1]^d:Az=Ay\}\) has a vertex \(z\) by Krein--Milman.
Lemma~\ref{lem:vertex-criterion} makes the free columns at \(z\) independent;
since \(A\) is singular, \(|F(z)|\le d-1\). Rounding at \(z\) gives
\[
\|A\xi-Ay\|_1\le
\sum_{i\in F(z)}(1-|z_i|)\le d-1,
\]
contrary to the hypothesis. Hence \(A\) is invertible. If some
\(|y_i|=1\), rounding \(y\) itself gives a sign vector \(\xi\) with
\(\|A(\xi-y)\|_1\le d-1\), again a contradiction. Thus
\(y\in(-1,1)^d\).

Put \(p=-Ay\). Then the origin is an interior point of the parallelepiped
\(A[-1,1]^d+p\). Fix \(r\in\{1,\dots,d\}\) and \(s\in\{-1,1\}\). Let
\(c=A^{-1}(s e_r)\), and for \(t\ge0\) set \(x(t)=y+tc\). Then
\[
Ax(t)+p=Ay+tA c+p=t s e_r.
\]
Let \(t_0>0\) be the first time at which a coordinate of \(x(t)\) reaches the
boundary of \([-1,1]^d\):
\[
t_0=\min_{i:\,c_i\ne0}\frac{\sgn(c_i)-y_i}{c_i}.
\]
Write \(x=x(t_0)\), and let \(i\) be a coordinate with
\(|x_i|=1\); put \(\eps=x_i\). Round \(x\) to a nearest sign vector \(\xi\),
chosen so that \(\xi_i=\eps\). By the hypothesis and \(\|a_k\|_1=1\),
\begin{equation}\label{eq:local-rigidity-key}
\begin{aligned}
d-1
&<\|A(\xi-y)\|_1
\le\|A(x-y)\|_1+\|A(\xi-x)\|_1\\
&=\|t_0s e_r\|_1+\left\|\sum_{k=1}^d(\xi_k-x_k)a_k\right\|_1
\le t_0+\sum_{k=1}^d|\xi_k-x_k|\\
&=t_0+\sum_{k=1}^d(1-|x_k|)
=t_0+(d-1)-\sum_{k\ne i}|x_k|\\
&\Longrightarrow\quad\sum_{k\ne i}|x_k|<t_0.
\end{aligned}
\end{equation}

Now set \(v=\eps a_i\) and \(z=p+v\). Since
\(Ax=\sum_k x_k a_k=\eps a_i+\sum_{k\ne i}x_k a_k\) and
\(Ax=t_0s e_r-p\),
\[
t_0s e_r-z=Ax-\eps a_i=\sum_{k\ne i}x_k a_k.
\]
By \eqref{eq:local-rigidity-key},
\(\|t_0s e_r-z\|_1<t_0\). On the other hand,
\[
\|t_0s e_r-z\|_1=|t_0s-z_r|+\sum_{k\ne r}|z_k|
\ge t_0-sz_r+\sum_{k\ne r}|z_k|.
\]
Hence \(sz_r>\sum_{k\ne r}|z_k|\), that is,
\(p+\eps a_i\in\mathcal C_{r,s}\).

The direction \(s e_r\) was arbitrary. Thus every cone
\(\mathcal C_{r,s}\) contains one of the points \(p+a_i\) or \(p-a_i\).
By Lemma~\ref{lem:cone-counting}, the dominant positions of the \(a_i\) are
pairwise distinct. Since \(\|a_i\|_1=1\), membership in a cone means that the
corresponding coordinate has absolute value larger than \(1/2\). This completes
the proof.
\end{proof}

The proof used only \(2d\) particular corners. Consequently, let \(A\) be
invertible with \(\|a_i\|_1=1\), and let \(0\) lie in the interior of
\(A[-1,1]^d+p\); equivalently, \(p=-Ay\) for some \(y\in(-1,1)^d\). If the
columns do not have pairwise distinct dominant positions, then one of the
following \(2d\) roundings satisfies \(\|A\xi+p\|_1\le d-1\): for each
\(s e_r\), move from \(y\) along \(A^{-1}(s e_r)\) to the first boundary and
round. Otherwise the rigidity argument would force distinct dominant
positions.

\section{Vertex exchanges and the parity contradiction}\label{sec:exchange}

We now describe exchanges between neighboring vertices, keeping assumption
\eqref{eq:H}. Let \(x\) be a vertex of \(K(V)\), and write
\(F(x)=\{i_1<\dots<i_d\}\). Put
\(a_j=v_{i_j}\), \(A=[a_1,\dots,a_d]\), and \(y_j=x_{i_j}\). Thus
\(A\) is invertible and \(\|y\|_1<1\) by the facts recalled above.

Choose a fixed index \(k\in Z(x)\), and put \(u=v_k\) and \(s=x_k\in\{-1,1\}\).
Let \(b=A^{-1}u\). Define
\begin{equation}\label{eq:move}
x_k(t)=s(1-t),\qquad
x_{i_j}(t)=y_j+tsb_j\quad(1\le j\le d),
\end{equation}
and leave all other coordinates unchanged. Since
\[
-su+\sum_{j=1}^d sb_j a_j=-su+sAb=0,
\]
we have \(Vx(t)=0\) for all \(t\). For \(t=0\) this is the vertex \(x\). We
move in the positive direction \(t\ge0\) and stop at the first time at which a
coordinate in \(F(x)\cup\{k\}\) reaches the boundary of the box, also allowing
the time \(t=2\); the coordinate \(k\) itself reaches the opposite boundary at
\(t=2\). Let \(t_0\) denote this first boundary time.

We first show that \(t_0<2\). If \(t_0=2\) and a free coordinate reaches the
boundary at \(t=2\), then \(x'=x(2)\) has fewer than \(d\) free coordinates.
Its free vectors form a subfamily of the independent \(a_1,\dots,a_d\), so
\(x'\) is a vertex by Lemma~\ref{lem:vertex-criterion}, contradicting the fact
that every vertex has \(d\) free coordinates. Hence no free coordinate is
fixed at \(t=2\). Then \(x'\) still has free coefficients \(y'=y+2sb\) and is
a vertex; from \(\|y\|_1,\|y'\|_1<1\) and
\(b=s(y'-y)/2\),
\[
1=\|u\|_1=\|Ab\|_1\le\|b\|_1\le\frac{\|y\|_1+\|y'\|_1}{2}<1,
\]
a contradiction. Thus \(t_0<2\).

At least one free coordinate reaches the boundary at \(t_0\). Let
\(Q=\{j:b_j\ne0\text{ and }y_j+t_0sb_j\in\{-1,1\}\}\). Since \(t_0<2\), the
coordinate \(k\) is free at \(x'\). The free coordinates at \(x'\) are
therefore \((F(x)\setminus Q)\cup\{k\}\). The corresponding vectors are a
subset of \(\{a_l:l\notin Q\}\cup\{u\}\). To see that this family is
independent, suppose
\[
c_0u+\sum_{l\notin Q}c_l a_l=0.
\]
Substituting \(u=\sum_{m=1}^d b_m a_m\) and using independence of the \(a_m\),
the coefficient of \(a_m\) for \(m\in Q\) is \(c_0b_m\), so \(c_0=0\); then
all \(c_l=0\). Hence the free vectors at \(x'\) are independent, and \(x'\) is
a vertex by Lemma~\ref{lem:vertex-criterion}. Every vertex has exactly \(d\)
free coordinates, so
\[
(d-|Q|)+1=d,
\]
which gives \(|Q|=1\). Thus exactly one old free vector has become fixed, the
new free family at \(x'\) is \(\{a_l:l\ne j\}\cup\{u\}\) for some \(j\), and
\(x'\) is the neighboring vertex obtained by the release.

Consequently every \(v_i\) has a dominant coordinate: apply
Lemma~\ref{lem:rigidity} when \(i\) is free at a vertex, and after releasing
\(i\) when it is fixed.

For each \(i\), choose \(\sigma_i\in\{-1,1\}\) so that
\(w_i=\sigma_i v_i\) has positive dominant coordinate; uniqueness makes this
well defined. Let
\(W=[w_1,\dots,w_n]\). The set of all signed sums is unchanged, since
\[
\sum_i\eps_i v_i=\sum_i(\eps_i\sigma_i)w_i.
\]
Hence \eqref{eq:H} remains true for \(W\). It suffices to prove the theorem for
the sign-normalized family \(w_1,\dots,w_n\). From now on we assume that every
\(w_i\) has its dominant coordinate positive, and we write \(v_i=w_i\) to avoid
extra notation.

\begin{lemma}[Invariance of the positive fixed count]\label{lem:P-invariant}
Let \(x\) and \(x'\) be the vertices related by the release construction above.
Define
\[
P(x)=\#\{i:x_i=1\}.
\]
Then \(P(x')=P(x)\).
\end{lemma}

\begin{proof}
Let \(\rho(a)\) denote the dominant position of a unit vector \(a\). The free
families at \(x\) and \(x'\) both have \(d\) vectors with pairwise distinct
dominant positions. Since they share \(\{a_l:l\ne j\}\), exactly \(a_j\) is
replaced by \(u\), and therefore
\[
\{\rho(a_1),\dots,\rho(a_d)\}=[d],\qquad \rho(u)=\rho(a_j).
\]
Relabel so that \(\rho(a_i)=i\); then \(\rho(u)=j\). In particular,
\(a_i\in\mathcal C_{i,+1}\) and \(u\in\mathcal C_{j,+1}\).

For \(t\in[0,1]\), put \(a_i(t)=(1-t)a_i+t a_{i,i}e_i\). Every matrix
\([a_1(t),\dots,a_d(t)]\) is strictly column diagonally dominant, hence
nonsingular, and the homotopies below preserve this property. Let
\(A[j\leftarrow u]\) denote the matrix obtained by replacing the \(j\)-th
column by \(u\). Homotoping all columns of \(A\) to diagonal form, and all
columns of \(A[j\leftarrow u]\) except the \(j\)-th, gives
\[
\det A=\prod_{i=1}^d a_{i,i}>0,\qquad
\det A[j\leftarrow u]=u_j\prod_{i\ne j}a_{i,i}>0.
\]
Cramer's rule therefore yields
\begin{equation}\label{eq:cramer}
b_j=(A^{-1}u)_j=\frac{\det A[j\leftarrow u]}{\det A}>0.
\end{equation}

Since \(x_{i_j}(t)=y_j+tsb_j\) and \(b_j>0\), the coordinate that reaches the
boundary at \(t_0\) has sign \(s\), so \(x'_{i_j}=s\). Also \(t_0<2\), hence
\(k\) is free at \(x'\), while all other fixed coordinates are unchanged.
Using indicator notation,
\[
P(x')-P(x)
=\mathbf{1}_{\{x'_{i_j}=1\}}-\mathbf{1}_{\{x_k=1\}}
=\mathbf{1}_{\{s=1\}}-\mathbf{1}_{\{s=1\}}=0.
\]
\end{proof}

Every edge of \(K(V)\) is a one-dimensional face. Passing from a vertex to an
incident edge keeps the \(d\) equations \(Vx=0\) active and drops exactly one
active box constraint, so it is the release move \eqref{eq:move}. Hence \(P\)
is constant on each edge by Lemma~\ref{lem:P-invariant}. Since the vertex graph
is connected \cite[Chapter~3]{Ziegler}, \(P\) is constant on all vertices.

Negation preserves convex combinations, so \(K(V)\) is centrally symmetric and
\(-x\) is a vertex.

The fixed coordinate set of \(-x\) is the same as that of \(x\), but all fixed
signs are reversed. Since \(x\) has \(n-d\) fixed coordinates,
\[
P(-x)=\#\{i:(-x)_i=1\}=\#\{i:x_i=-1\}=(n-d)-P(x).
\]
On the other hand, \(P(x)=P(-x)\) because \(P\) is constant on vertices.
Therefore \(n-d=2P(x)\) is even. But \(n\not\equiv d\pmod2\) means that
\(n-d\) is odd. This contradiction proves that \eqref{eq:H} is impossible, and
hence establishes \eqref{eq:main}. Undoing the sign normalizations, if any,
gives the stated signs for the original vectors.

\begin{remark}[The target must be the origin]
The same argument does not apply to a nonzero fiber
\(\{x\in[-1,1]^n:Vx=c\}\) with \(c\ne0\). The map \(x\mapsto-x\) sends that
fiber to the fiber with target \(-c\), not back to itself, so one does not get
the equation \(P(x)=P(-x)\) inside one polytope. This is consistent with the
example in the introduction.
\end{remark}

\section{Deterministic constructions}\label{sec:algorithm}

The proof above is existential. We give a deterministic \(O(nd^2)\)
construction, matching classical rounding. Section~\ref{subsec:blocked}
reduces this below \(nd^2\) by blocking and fast matrix multiplication.

Choose \(\sigma_i\in\{-1,1\}\) so that \(w_i=\sigma_i v_i\) has a positive
dominant coordinate when one exists, and set \(\sigma_i=1\) otherwise. For
\(W=[w_1,\dots,w_n]\), signs \(\zeta\) for \(W\) give signs
\(\eps_i=\sigma_i\zeta_i\) for \(V\).

Put \(M=\begin{bmatrix}W\\ \one^T\end{bmatrix}\), where
\(\one=(1,\dots,1)^T\), and consider
\begin{equation}\label{eq:K0}
K_0=\{z\in[-1,1]^n:Wz=0,\ \one^Tz=0\}.
\end{equation}
This is the zero fiber of \(M\); by Section~\ref{sec:rounding}, it has a
vertex.

\begin{proposition}[From a vertex of \(K_0\) to a good rounding]\label{prop:endpoints}
Let \(z\) be a vertex of \(K_0\). Either direct rounding at \(z\) gives
\(\|W\zeta\|_1\le d-1\), or one endpoint of the edge of \(K(W)\) through \(z\)
satisfies the local-search hypothesis.
\end{proposition}

\begin{proof}
Let \(F=F(z)\) and \(f=|F|\). Since the free columns of \(M\) are independent,
\(f\le d+1\). The rounding estimate \eqref{eq:rounding}, applied to \(W\), gives
\begin{equation}\label{eq:K0-round}
\|W\zeta\|_1\le f-\|z_F\|_1.
\end{equation}

If \(f\le d-1\), the right side is at most \(d-1\). If \(f=d\), the \(n-d\)
fixed signs form a nonzero odd integer, hence have absolute value at least
\(1\). Since \(\one^Tz=0\),
\[
\|z_F\|_1\ge\left|\sum_{i\in F}z_i\right|
=\left|\sum_{i\notin F}z_i\right|\ge1.
\]
Thus \eqref{eq:K0-round} is at most \(d-1\).

It remains to consider \(f=d+1\). Let \(W_F\) be the \(d\times(d+1)\) matrix
formed by the free columns. Its kernel is one-dimensional: for nonzero
\(b\in\ker W_F\), independence of the free columns of \(M\) gives
\(\one^Tb\ne0\). If \(h\in\ker W_F\), then
\(h'=h-(\one^Th/\one^Tb)b\) satisfies \(W_Fh'=0\) and \(\one^Th'=0\); hence
\(M_Fh'=0\), so \(h'=0\). Thus \(h\) is a multiple of \(b\).

Moving from \(z_F+tb\) in both directions to the first boundary gives endpoints
\(x^\pm\in K(W)\), with parameters \(t^\pm\). Since \(\one^Tz=0\),
\[
\one^Tx^+=t^+\one^Tb,\qquad \one^Tx^-=t^-\one^Tb,
\]
so \((\one^Tx^-)(\one^Tx^+)<0\).

At each endpoint a free coordinate \(i\) has \(b_i\ne0\). Deleting it leaves
\(d\) independent columns, since otherwise there is a kernel vector with
\(i\)-th coordinate zero, impossible as all kernel vectors are multiples of
\(b\). Thus both endpoints are vertices of \(K(W)\).

If at an endpoint the number of free coordinates is less than \(d\), or the
\(\ell_1\)-norm of the free coefficients is at least \(1\), direct rounding at
that endpoint gives the bound. Otherwise both endpoints have exactly \(d\) free
coordinates and free coefficient vectors \(y^\pm\) with
\(\|y^\pm\|_1<1\). Let \(T^\pm\) be the sum of the \(n-d\) fixed signs at the
two endpoints. Since \(n-d\) is odd, each \(T^\pm\) is a nonzero odd integer,
so \(|T^\pm|\ge1\). Write \(S^\pm\) for the sum of the free coefficients at
the two endpoints. Then
\[
\one^Tx^\pm=T^\pm+S^\pm,\qquad |S^\pm|\le\|y^\pm\|_1<1.
\]
Thus \(\one^Tx^\pm\) has the same sign as \(T^\pm\). Since the two numbers
\(\one^Tx^\pm\) have opposite signs, \(T^-T^+<0\).

They cannot fix the same free coordinate: if both fixed \(i\), then this
coordinate differs by \(2\) between \(x^-\) and \(x^+\). With \(A'\) denoting
the common other free vectors, subtracting \(Wx^\pm=0\) gives
\[
\pm2w_i+A'(y^+-y^-)=0.
\]
Taking \(\ell_1\)-norms and using \(\|w_i\|_1=1\),
\[
2\le 2\|w_i\|_1=\|A'(y^+-y^-)\|_1\le\|y^+\|_1+\|y^-\|_1<2,
\]
a contradiction. Hence the endpoints exchange two free coordinates, say
\(i\) and \(j\).

Finally, the two endpoint families cannot both be rigid. If they were, the
difference \(x^+-x^-\) lies in the one-dimensional kernel and has the sign of
\(-x^-_i\), so releasing \(i\) from \(x^-\) gives \(x^+\). Lemma~\ref{lem:P-invariant}
then makes the entering and leaving vectors share a dominant position and gives
a positive coefficient; hence the fixed signs, and \(T^\pm\), are unchanged,
contradicting \(T^-T^+<0\). Thus one endpoint is not rigid, and the local-search
observation gives a good candidate.
\end{proof}

To find a vertex of \(K_0\), perform Gaussian elimination on \(M\) and obtain
independent rows \(U\in\R^{r\times n}\), \(r\le d+1\), with the same nullspace.
Since \(\rank U=r\), choose \(r\) independent columns; they form an invertible
\(B\), indexed by \(J\). Set \(C=B^{-1}\) and \(x=0\).

Process \(k\notin J\) in order. Coordinates in \(J\) are allowed to lie in
\([-1,1]\), processed coordinates outside \(J\) are fixed at \(\pm1\), and
unprocessed ones are \(0\). For the next \(k\), let \(u_k\) be the \(k\)-th
column of \(U\) and set
\[
\alpha=Cu_k,\qquad b_J=-\alpha,\qquad b_k=1,
\]
with all other coordinates zero. Then \(Ub=0\). Move \(x+tb\) to the first
boundary,
\[
t_0=\min\left\{1,\ \min_{j\in J:\,b_j\ne0}\frac{\sgn(b_j)-x_j}{b_j}\right\},
\]
with empty minimum \(+\infty\). If \(x_k\) reaches \(1\) first, fix it at
\(1\); otherwise a coordinate \(j\in J\) reaches its boundary, is fixed there,
and \(u_k\) replaces \(B_j\). Break ties by a fixed order.

Since \(\alpha_j\ne0\), the replacement stays invertible: a dependence would
have coefficient \(c_j\alpha_j\) at the exiting column, forcing \(c_j=0\) and
then all coefficients zero.

The inverse is updated in place. If \(B_j\) is replaced, then
\begin{equation}\label{eq:inverse-update}
C'=C-\frac{\alpha-e_j}{\alpha_j}\,C_{j,:},
\end{equation}
where \(C_{j,:}\) is row \(j\) of \(C\). Indeed, for unchanged columns \(B_i\),
\(CB_i=e_i\) and \(C_{j,:}B_i=0\); for \(u_k=B\alpha\),
\(Cu_k=\alpha\) and \(C_{j,:}u_k=\alpha_j\). Hence \(C'B'=I\).

Each update costs \(O(d^2)\). Index \(k\) is processed even when \(t_0=0\), so
there are exactly \(n-r\) rounds and no cycling. At the end all nonbasis
coordinates are fixed and the free columns lie in \(B\);
Lemma~\ref{lem:vertex-criterion} makes \(x\) a vertex of \(K_0\).

For an endpoint selected by Proposition~\ref{prop:endpoints}, let \(A\) be its
free matrix, \(y\) its free vector, and \(p=-Ay\). For each direction
\(s e_r\), set \(b=sA^{-1}e_r\) and move to the first boundary:
\[
x=y+tb,\qquad t=\min_{i:b_i\ne0}\frac{\sgn(b_i)-y_i}{b_i}.
\]
Round \(x\) to \(\xi\), with \(0\mapsto1\), and test
\(\|A\xi+p\|_1=\|A(\xi-y)\|_1\). The local-search observation guarantees success
among the at most \(2d\) candidates. Since \(A^{-1}\) costs \(O(d^3)\) and each
test \(O(d^2)\), this step costs \(O(d^3)\).

The complete algorithm is therefore:
\begin{enumerate}
\item If \(n<d\), take all signs equal to \(+1\). Otherwise adjust the vectors
to obtain \(W\) and the signs \(\sigma_i\).
\item Find a vertex \(z\) of \(K_0\) by the basis procedure above.
\item Apply Proposition~\ref{prop:endpoints}: round directly, use one of the two
endpoints, or perform the local search with at most \(2d\) candidates.
\item Undo the sign adjustment and output \(\eps_i=\sigma_i\zeta_i\).
\end{enumerate}

The initial elimination and inverse cost \(O(nd^2)\), and the main loop has
\(n\) rounds of cost \(O(d^2)\). Since the endpoint step costs
\(O(d^3)\le O(nd^2)\) for \(n>d\), the total is
\begin{equation}\label{eq:complexity}
O(nd^2),
\end{equation}
with \(O(nd)\) storage.

\begin{theorem}\label{thm:algorithm-direct}
Under the hypotheses of Theorem~\ref{thm:main}, the signs in \eqref{eq:main}
can be constructed deterministically in \(O(nd^2)\) exact arithmetic operations
and comparisons.
\end{theorem}

\subsection{Blocked basis maintenance}\label{subsec:blocked}

Blocking removes the \(O(d^2)\) per-round inverse update. Let \(\omega\) be a
matrix-multiplication exponent. Strassen gives \(\omega<3\) \cite{Strassen};
the bound \(\omega\le9/4\) over \(\mathbb C\) is due to \cite{OpenAI26} and holds
over \(\R\) by embedding, and \cite{Lu26} extends it to arbitrary
characteristic.

Rank-profile decomposition computes the full-row-rank matrix \(U\in\R^{r\times n}\)
with \(r\le d+1\), independent columns \(J\), and \(C=U_J^{-1}\) in
\(O(nd^{\omega-1}+d^\omega)\) field operations \cite{BunchHopcroft,DPS16}. For a
pending \(k\), the query and update are
\begin{equation}\label{eq:block-query}
\alpha=CU_k,\qquad
C'=C-\frac{\alpha-e_j}{\alpha_j}\,e_j^TC,
\end{equation}
where \(j\) is the leaving index; the direct update costs \(O(r^2)\).

The update is left-multiplicative. For a block starting with inverse \(C_0\),
after \(p\) pivots,
\begin{equation}\label{eq:block-low-rank}
C_p=C_0-\sum_{t=1}^p v_tw_t^T,
\end{equation}
where, if the position \(j_t\) leaves at step \(t\),
\begin{equation}\label{eq:block-pairs}
v_t=\frac{\alpha^{(t)}-e_{j_t}}{\alpha^{(t)}_{j_t}},\qquad
w_t^T=e_{j_t}^TC_{t-1}.
\end{equation}
Thus a query need not form \(C_p\):
\begin{equation}\label{eq:block-eval}
C_pU_k=C_0U_k-\sum_{t=1}^p v_t\,(w_t\cdot U_k).
\end{equation}
Moreover, \eqref{eq:block-low-rank} gives
\begin{equation}\label{eq:block-w}
w_t=C_0^Te_{j_t}-\sum_{s<t}(v_s)_{j_t}w_s.
\end{equation}

\begin{lemma}[Block simulation]\label{lem:block-simulation}
Let \(B\le r\). Given \(C_0\) and the columns in a block of size \(B\), one can
simulate all queries and pivots of that block, and recover the exact inverse at
its end, in \(O(r^\omega+B^2r)\) field operations. With the same tie-breaking
rule, the simulated pivot sequence is identical to the direct sequence.
\end{lemma}

\begin{proof}
One padded \(r\times r\) product computes all \(C_0U_k\) in \(O(r^\omega)\).
Each query \eqref{eq:block-eval} and each \(w_t\) in \eqref{eq:block-w} costs
\(O(Br)\); with at most \(B\) of each, these total \(O(B^2r)\). A second padded
product collapses \eqref{eq:block-low-rank} in \(O(r^\omega)\). Since
\eqref{eq:block-eval} and \eqref{eq:block-w} are exact, all comparisons and
pivots match the direct walk.
\end{proof}

Use consecutive blocks of size
\(B=\min\{r,\lceil (d+1)^{(\omega-1)/2}\rceil\}\); smaller \(r\) only improves
the bound. By Lemma~\ref{lem:block-simulation}, the basis walk costs
\begin{equation}\label{eq:block-time}
O\!\left(\frac nB(d+1)^\omega+nB(d+1)\right).
\end{equation}
The terms balance at \(B\), both with exponent \((\omega+1)/2\) in \(d\).
Hence the vertex of \(K_0\) is found in
\begin{equation}\label{eq:vertex-time}
O\!\left(nd^{(\omega+1)/2}+d^\omega\right)
\end{equation}
operations, and the rank-profile cost is dominated since \(\omega<3\).

Batch the endpoint test. For the endpoint matrix \(A\) and free vector \(y\),
set \(p=-Ay\) and compute \(A^{-1}\) in \(O(d^\omega)\). For \(se_r\), let
\(z=y+tb\) be the first boundary point; the rounded vector \(\xi\) satisfies
\begin{equation}\label{eq:batched-endpoint}
A\xi+p=A(\xi-z)+tse_r.
\end{equation}
Stack the \(2d\) vectors \(\xi-z\) into two padded \(d\times d\) blocks. Two fast
multiplications by \(A\) compute all \(2d\) terms in \(O(d^\omega)\); shifts and
\(\ell_1\)-norms add \(O(d^2)\). Thus the endpoint and local-search step costs
\(O(n+d^\omega)\).

\begin{theorem}\label{thm:algorithm-fast}
Assume \(n>d\). The signs in \eqref{eq:main} can be constructed
deterministically in
\begin{equation}\label{eq:fast-complexity}
O\!\left(nd^{(\omega+1)/2}+d^\omega\right)
\end{equation}
exact arithmetic operations and comparisons. With the bound \(\omega\le9/4\),
this is \(O(nd^{13/8})\), which is strictly below the direct \(O(nd^2)\) bound.
\end{theorem}

\begin{proof}
Correctness follows from the direct construction: Lemma~\ref{lem:block-simulation}
reproduces the pivot sequence and \eqref{eq:batched-endpoint} evaluates the same
candidates. The vertex walk, endpoint step, and sign adjustment cost
\eqref{eq:vertex-time}, \(O(n+d^\omega)\), and \(O(nd)\). Since \(n>d\) and
\(\omega<3\), the \(d^\omega\) term is absorbed by \(nd^{(\omega+1)/2}\), giving
\eqref{eq:fast-complexity}. Finally, \((\omega+1)/2\le13/8<2\) for
\(\omega\le9/4\); if \(n<d\), all-\(+1\) costs \(O(n)\).
\end{proof}

\begin{remark}
The gain is asymptotic: fast-matrix-multiplication constants are large, and no
bit-complexity is claimed; all bounds use exact real arithmetic.
\end{remark}

\appendix

\section{Height interpretation}\label{app:height}

This appendix gives a heuristic interpretation. Attach heights \(h_i\in\R\) to
the vectors \(v_i\), and consider the zonotope
\[
Z_h=\sum_{i=1}^n[-(v_i,h_i),(v_i,h_i)].
\]
Its points with vanishing first \(d\) coordinates are exactly \((0,h\cdot x)\)
with \(x\in K(V)\), so the lowest height is
\[
H(h)=\min_{x\in K(V)}h\cdot x.
\]
For generic \(h\), the minimizer is a unique vertex and changes along an edge
as \(h\) varies. The exchange is thus a switch of the lowest representation.

After the sign adjustment, label a directly unroundable vertex \(x\) with
pairwise distinct dominant positions by
\[
\sgn\bigl(2P(x)-(n-d)\bigr).
\]
This label is nonzero, constant on adjacent vertices of this type by
Lemma~\ref{lem:P-invariant}, and negated at \(-x\). Along a height path from
\(x\) to \(-x\) it must change; at a change point the vertex is directly
roundable or covered by the local-search observation.

\renewcommand{\bibfont}{\small\raggedright}
\setlength{\bibsep}{0pt}
\bibliographystyle{plain}
\bibliography{references}

\end{document}